\documentclass[a4paper,10pt]{amsart}

\usepackage{verbatim}
\usepackage{enumerate}
\usepackage{amsmath, amsthm, amssymb}
\usepackage{amsfonts,bm}
\usepackage{hyperref}

\newtheorem{thm}{Theorem}[section]
\newtheorem{lem}[thm]{Lemma}
\newtheorem{defn}[thm]{Definition}
\newtheorem{prop}[thm]{Proposition}
\newtheorem{cor}[thm]{Corollary}

\theoremstyle{remark}
\newtheorem{ex}[thm]{Example}
\newtheorem{rem}[thm]{Remark}

\newcommand{\beq}{\begin{equation}}
\newcommand{\eeq}{\end{equation}}
\newcommand{\bal}{\begin{aligned}}
\newcommand{\eal}{\end{aligned}}

\newcommand{\ben} {\begin{enumerate}}
\newcommand{\beni} {\begin{enumerate}[(i)]}
\newcommand{\een} {\end{enumerate}}
\newcommand{\beqw} {\begin{equation*}}
\newcommand{\eeqw} {\end{equation*}}

\numberwithin{equation}{section}

\renewcommand{\a}{\alpha}

\renewcommand{\b}{\beta}
\renewcommand{\d}{\delta}
\newcommand{\g}{\gamma}

\renewcommand{\l}{\lambda}

\renewcommand{\O}{\Omega}

\newcommand{\s}{\sigma}

\newcommand{\E}{\mathbb{E}}
\newcommand{\R}{\mathbb{R}}

\newcommand{\C}{\mathbb{C}}

\newcommand{\F}{\mathcal{F}}
\renewcommand{\P}{\mathbb{P}}
\renewcommand{\L}{\mathcal{L}}
\newcommand{\B}{\mathcal{B}}
\newcommand{\dD}{\mathbb{D}}
\newcommand{\eE}{\mathcal{E}}

\newcommand{\one}{\mathbf{1}}

\newcommand{\bpr}{\begin{proof}}
\newcommand{\epr}{\end{proof}}
\newcommand{\blem}{\begin{lem}}
\newcommand{\elem}{\end{lem}}
\newcommand{\bdefn}{\begin{defn}}
\newcommand{\edefn}{\end{defn}}
\newcommand{\bthm}{\begin{thm}}
\newcommand{\ethm}{\end{thm}}
\newcommand{\bprop}{\begin{prop}}
\newcommand{\eprop}{\end{prop}}
\newcommand{\bcor}{\begin{cor}}
\newcommand{\ecor}{\end{cor}}
\newcommand{\bex}{\begin{ex}}
\newcommand{\eex}{\end{ex}}
\newcommand{\brem}{\begin{rem}}
\newcommand{\erem}{\end{rem}}

 \bmdefine{\i}{i}

\renewcommand{\H}{\mathcal{H}}
\newcommand{\calS}{\mathcal{S}}
\newcommand{\calP}{\mathcal{P}}
\newcommand{\Dom}{{\mathsf D}}

\newcommand{\ip}[1]{\langle {#1}\rangle}

\newcommand{\lin}{\mathrm{lin}}

\newcommand{\ot}{\otimes}
\newcommand{\hot}{\widehat\otimes}
\newcommand{\sot}{{\hbox{\tiny\textcircled{s}}}}

\author{Jan Maas}
\address{Delft Institute of Applied Mathematics\\
Delft University of Technology\\ P.O. Box 5031\\ 2600 GA Delft\\The
Netherlands } \email{J.Maas@tudelft.nl}

\thanks{The author acknowledges support by the `VIDI subsidie' 639.032.201
of the Netherlands Organisation for Scientific Research (NWO) and the
ARC Discovery Grant DP0558539}



\title[Malliavin calculus and decoupling inequalities]{Malliavin calculus and decoupling inequalities in Banach spaces}

\begin{document}

\begin{abstract} We develop a theory of Malliavin calculus for
Banach space valued random variables. Using radonifying operators
instead of symmetric tensor products we extend the Wiener-It\^o
isometry to Banach spaces. In the white noise case we obtain two
sided $L^p$-estimates for multiple stochastic integrals in arbitrary
Banach spaces. It is shown that the Malliavin derivative is bounded
on vector-valued Wiener-It\^o chaoses. Our main tools are decoupling
inequalities for vector-valued random variables. In the opposite
direction we use Meyer's inequalities to give a new proof of a
decoupling result for Gaussian chaoses in UMD Banach spaces.
\end{abstract}


\maketitle

\section{Introduction}

The theory of Malliavin calculus \cite{IW89,Nu06} has been developed
in the seventies by Malliavin \cite{Ma78}, who used it to give a
probabilistic proof of H\"ormander's ``sums of squares''-theorem. The
Malliavin calculus generalises in a natural way to Hilbert space
valued random variables. We refer to  \cite{CT06} for a recent
account of this infinite dimensional setting with applications to
stochastic (partial) differential equations.

In recent years many Hilbert space results in stochastic (and
harmonic) analysis have been transferred to a Banach space setting
\cite{hytonen:lps,KWsf}. Of particular relevance for this work is the
theory of stochastic integration in Banach spaces developed by van
Neerven, Veraar and Weis \cite{vNVW07,vNW05}. Motivated by these
developments we construct in this paper a theory of Malliavin
calculus for random variables taking values in a Banach space.

Vector-valued Malliavin calculus has been consider by several authors
\cite{MN93,MWZ05,Sh94}. The main focus in this work is on the
interplay between Malliavin calculus and decoupling inequalities. On
the one hand, decoupling inequalities are our main tools in the proof
of Theorems \ref{thm:WienerIto}, \ref{thm:multstochint} and
\ref{thm:Dchaosm}. In the opposite direction, we apply the theory
developed in this paper to give a new proof of a known decoupling
result in Theorem \ref{thm:decouplingUMD}.

In a follow-up paper with van Neerven \cite{MvNCO} the vector-valued
Malliavin calculus is used to construct a Skorokhod integral in UMD
spaces which extends the stochastic integral from \cite{vNW07}. This
is used to obtain a Clark-Ocone representation formula in UMD spaces.

It has been proved by Pisier \cite{Pi88} that the fundamental Meyer
inequalities remain valid if the Banach space is a UMD space,
provided that the norm of the derivative is taken in the appropriate
space. These spaces turn out to be spaces of so-called
$\g$-radonifying operators, which have been used to transfer
classical Hilbert space results to a more general Banach space
setting in various recent works.

Firstly, in the work of Kalton and Weis on $H^\infty$-functional
calculus \cite{KW04} $\g$-radonifying operators appear as
generalisations of classical square functions from harmonic analysis.
Secondly, $\g$-radonifying operators are used in \cite{vNVW07,vNW05}
to obtain two-sided estimates for moments of vector-valued stochastic
integrals, and provide a generalisation of the classical
It\^o-isometry.

Let us describe some of the main results in this paper. For details
we refer to later sections. Let $(\O,\F,\P)$ be a probability space,
let $H$ be a Hilbert space, and let $E$ be a Banach space. We
consider an isometry $W: H \to L^2(\O)$ onto a closed subspace
consisting of Gaussian random variables, and assume that $\F$ is the
$\s$-field generated by $\{W(h) : h \in H\}.$ The classical
Wiener-It\^o decomposition says that $L^2(\O,\F,\P)$ admits an
orthogonal decomposition into Gaussian chaoses $L^2(\O,\F,\P) =
\bigoplus_{m\geq0}\H_m.$ Moreover, there exist canonical isometries
 $\Phi_m$
from the symmetric Hilbert space tensor powers $H^{\sot m}$ onto
$\H_m.$

We show in Theorem \ref{thm:WienerIto} that this result admits a
natural Banach space valued generalisation. For this purpose we
consider the space of symmetric $\g$-radonifying operators $\g^{\sot
m}(H,E)$ (cf. Section \ref{sec:chaos}), which turns out to be the
natural vector-valued analogue of the symmetric Hilbert space tensor
powers. We prove that $\Phi_m$ extends to an $L^p$-isomorphism
between $\g^{\sot m}(H,E)$ and the vector-valued Gaussian chaos
$\H_m(E)$ for $1\leq p<\infty,$
 $$ \| (\Phi_m \ot I) T\|_{L^p(\O;E)}
   \eqsim_{m,p} \|T\|_{\g^{\sot m}(H,E)},
 \qquad T \in \g^{\sot m}(H,E).$$

In Section \ref{subsec:whitenoise} we consider the particular case
where $H = L^2(M,\mu)$ for some $\s$-finite measure space $(M,\mu).$
Theorem \ref{thm:multstochint} shows that the Wiener-It\^o
isomorphism between $\g^{\sot m}(L^2(M,\mu),E)$ and $\H_m(E)$ is
given by a multiple stochastic integral $I_m$ for Banach space valued
functions. This result gives two-sided bounds for $L^p$-norms of
multiple stochastic integrals, for $1\leq p<\infty,$
$$ \|I_m F\|_{L^p(\O;E)}
        \eqsim_{m,p} \|  F \|_{\g^{\sot m}(L^2(M),E)},\qquad
       F \in  \g^{\sot m}(L^2(M,\mu),E),$$
thereby generalising the (single) Banach space valued stochastic
integral of \cite{vNW05}.

The proofs of both results rely on (different) decoupling
inequalities. The idea to use decoupling in the study of multiple
stochastic integrals is not new. In fact, applications to multiple
stochastic integration appear already in the pioneering work on
decoupling by McConnell and Taqqu \cite{McCT86,McCT87}, Kwapie{\'n}
\cite{Kw87} and others. The decoupling results that we will use, as
well as a some preliminaries on $\g$-radonifying operators, can be
found in Section \ref{sec:prel}.

In Section \ref{sec:derivative} we consider the Banach space valued
Malliavin derivative $D,$ which for $1 \leq p< \infty$ acts as a
closed operator
 \begin{align*}
  D : \dD^{1,p}(E) \subset L^p(\O;E) \to L^p(\O;\g(H,E)).
 \end{align*}
The main result in this section (Theorem \ref{thm:Dchaosm}) asserts
that the restriction of the Malliavin derivative to each chaos is an
$L^p$-isomorphism for $1\leq p<\infty,$ $$\|DF\|_{L^p(\O;\g(H,E))}
\eqsim_{p,m}
  \|F\|_{L^p(\O;E)},\qquad F\in \H_m(E),
  $$
a fact which is by no means obvious for general Banach spaces.
The use of decoupling in this context appears to be new. In UMD
spaces this result is an easy consequence of Meyer's inequalities.

These inequalities are considered in more detail in Section
\ref{sec:OU}. We discuss several of its consequences and obtain a
version of Meyer's multiplier theorem in UMD spaces. We return to
decoupling in Theorem \ref{thm:decouplingUMD} where we give a new
proof of a known decoupling result for Gaussian chaoses in UMD
spaces based on
Meyer's inequalities. %

\medskip

{\em Acknowledgment} -- Part of this work was done during a half-year
stay of the author at the University of New South Wales in Sydney. He
thanks Ben Goldys for his kind hospitality.

\section{Preliminaries} \label{sec:prel}

\subsection{Decoupling} \label{subsec:decoupling}

Decoupling inequalities go back to the work of McConnell and Taqqu
\cite{McCT86,McCT87}, Kwapie{\'n} \cite{Kw87}, Arcones and Gin\'{e}
\cite{AG93}, and de la Pe\~na and Montgomery-Smith \cite{dlPMS95}
among others. We refer to the monographs \cite{KW92,dlPG99} for
extensive information on this topic.

First we introduce some notation which will be used throughout the
paper. For $j \geq 1$ and a finite sequence $\i = (i_1,  \ldots,
i_n)$ with values in $\{1,2,\ldots\}$ we set $$j(\i) = \#\{ i_k : 1
\leq k \leq n, i_k = j\}, \quad |\i| = n, \quad |\i|_\infty :=
\max_{1\leq k\leq n} i_k, \quad \i! = \prod_{k=1}^n k(\i)!.$$
Let $(\g_n)_{n \geq 1}$ be a Gaussian sequence, i.e. a sequence of
independent standard Gaussian random variables on a (sufficiently
rich) probability space $(\O,\F,\P),$ and let $(\g_n^{(k)})_{n \geq
1}$ be independent copies  for each $k \geq 1.$
The Hermite polynomials $H_m$ satisfy $H_0(x) = 1,$ $H_1(x) = x$ and
the recurrence relation
$$(m+1)H_{m+1}(x) = x H_m(x) - H_{m-1}(x), \quad m \geq 1.$$ We
set
$$\Psi_\i = (\i!)^{1/2} \prod_{j\geq 1}H_{j(\i)}(\g_j).$$

The next theorem states two well-known decoupling results which were
obtained in \cite{McCT87,Kw87,AG93}. A general result containing
both parts of the next theorem is due to Gin\'{e} \cite[Theorem
4.2.7]{dlPG99}.

We use the Vinogradov notation $\eqsim, \lesssim, \gtrsim$ to denote
estimates with universal constants. If the constants depend on
additional parameters, we include them as subscripts.

 \bthm \label{thm:decouplingD}
Let $E$ be a Banach space, let $m,n \geq 1,$ and suppose that we are
in one of the following two situations:
 \begin{enumerate}
    \item (symmetric case)
Let $(x_\i)_{|\i| = m} \subset E$ satisfy $x_\i = x_{\i'}$ whenever
$\i'$ is a permutation of $\i,$
 and set $$F := \sum_{|\i| =
m,|\i|_\infty \leq n} (\i!/m!)^{1/2}\Psi_\i x_{\i}.$$
   \item (tetrahedral case)
Let $(x_\i)_{|\i| = m} \subset E$ satisfy $x_\i = 0$ whenever $j(\i)
>1$ for some $j\geq 1,$ and set $$F := \sum_{|\i|= m,|\i|_\infty
\leq n}\g_{i_1} \cdot\ldots\cdot \g_{i_m} x_\i .$$
\end{enumerate}
In both cases we put $$\widetilde F := \sum_{|\i|= m,|\i|_\infty \leq
n}
     \g_{i_1}^{(1)} \cdot\ldots\cdot \g_{i_m}^{(m)}x_\i.$$
Then there exists a constant $C_m \geq 1$ depending only on $m,$ such
that for all $t>0$ we have
 \begin{align*}
 \frac{1}{C_m} \P(\|\widetilde F\|_E > C_m t )
 \leq \P(\| F\|_E >  t )
 \leq C_m \P(\|\widetilde F\|_E > \frac{t}{C_m} ).
 \end{align*}
Consequently, for $1 \leq p < \infty$ we have
 \begin{align*}
 \|F\|_{L^p(\O;E)} \eqsim_{p,m} \|\widetilde F\|_{L^p(\O;E)}.
 \end{align*}
 \ethm

\brem The requirement that $|\i|_\infty \leq n$ is chosen for
convenience, to ensure that we are dealing with finite sums
exclusively. Note however that the constants in all of our estimates
do not depend on $n.$ \erem

\subsection{Spaces of $\g$-radonifying operators}
\label{subsec:gammanorms}

In this section we will review some well-known results about
$\g$-radonifying operators. For more information we refer to
\cite{Bo,KWsf}. Let $H$ be a real separable Hilbert space with
orthonormal basis $(u_n)_{n \geq 1},$ and let $E$ be a real Banach
space. Let $(\g_n)_{n \geq 1}$ be a Gaussian sequence.

An operator $T \in \L(H,E)$ is said to be $\g$-radonifying if the sum
$\sum_{n = 1}^\infty \g_n T u_n$ converges in $L^2(\O;E).$ The
convergence and the $L^2$-norm of this sum do not depend on the
choice of the orthonormal basis and the Gaussian sequence. The space
$\g(H,E)$ consisting of all $\g$-radonifying operators in $\L(H,E)$
is a Banach space endowed with the norm
 \begin{align*}
 \|T\|_{\g(H,E)} := \bigg(\E \Big\| \sum_{n = 1}^\infty \g_n T u_n
 \Big\|_E^2\bigg)^{1/2}.
 \end{align*}
Obviously all  rank-$1$ operators
 $$h \ot x : h' \mapsto [h',h]\cdot x, \quad h,h' \in H, x \in E,$$
are contained in $\g(H,E)$ and one easily sees that they span a
dense subspace of $\g(H,E).$

An important role in this work will be played by spaces of the form
$\g^m(H,E),$ which we define inductively by
$$\g^1(H,E) := \g(H,E), \quad \g^{m+1}(H,E) := \g(H, \g^{m}(H,E)),
  \quad m \geq 1.$$
To improve readability we will write $T(h,h')$ instead of $(Th)(h')$
if $T \in \g^2(H,E).$ Furthermore we will write $ (h \ot h') \ot x$
to denote the operator $h \ot (h' \ot x) \in \g^2(H,E).$ Similar
remarks apply when $m>2.$
For future use we record that for operators of the form
 \begin{align} \label{eq:elementaryoperator}
 T = \sum_ {|\i| = m, |\i|_\infty \leq n}
 (u_{ i_1} \ot \cdots\ot u_{i_m}) \ot x_{\i}, \quad x_\i \in E,
 \end{align}
the norm in $\g^m(H,E)$ is given by
 \begin{align}  \label{eq:elementaryoperatornormgm}
 \|T\|_{\g^m(H,E)}^2
  = \E \Big\| \sum_{|\i| = m, |\i|_\infty \leq n}
 \g_{i_1}^{(1)} \cdot \ldots\cdot \g_{i_m}^{(m)}
   x_{\i}  \Big\|_E^2,
 \end{align}
where we use the multi-index notation from Section
\ref{subsec:decoupling}.

If $K$ is a Hilbert space then $\g^m(H,K)$ is canonically isometric
to the Hilbert space tensor product $H^{\hot m} \hot K.$ It has been
shown in \cite{KWsf} (see also \cite{vNW07}) that $\g^m(H,E)$ is
isomorphic to $\g(H^{\hot m},E)$ for all $m\geq 1$ if and only if the
Banach space $E$ has Pisier's property $(\a)$ \cite{Pi78}.

It is well known (see \cite[Section 5]{KWsf}) that the pairing
 \begin{align*}
[T,S]_{\g} :=  \text{tr\;}(T^*S),
  \quad  T \in \g(H,E), S\in \g(H,E^*),
 \end{align*}
defines a duality between $\g(H,E)$ and $\g(H,E^*),$ which allows us
to identify $\g(H,E^*)$ with a weak$^*$-dense subspace of the dual
space $\g(H,E)^*.$ It has been proved by Pisier \cite{Pis89} that the
Banach spaces $\g(H,E)^*$ and $\g(H,E^*)$ are isomorphic if $E$ is
$K$-convex. The notion of $K$-convexity and its relevance for
vector-valued Malliavin calculus will be discussed in Remark
\ref{rem:Kconvex}.
It is not difficult to check that $$[T,S]_{\g}= \sum_{j = 1}^\infty
\ip{Tu_j,Su_j} \quad \text{ and }\quad
 [T,S]_\g \leq  \|T\|_{\g(H,E)}
    \|S\|_{\g(H,E^*)}.$$

Let us now consider the important special case that $H = L^2(M,\mu)$
for some $\sigma$-finite measure space $(M,\mu).$ A strongly
measurable function $\phi:M^m \to E$ is said to be weakly-$L^2$ if
$\ip{\phi,x^*} \in L^2(M^m)$ for all $x^* \in E^*.$ We say that such
a function represents an operator $T_\phi \in \g^m(L^2(M),E)$ if
for all $f_1, \ldots, f_m \in L^2(M)$ and for all $x^* \in E^*$
we have
 \begin{align*}
 \ip{T_\phi (f_1, \ldots, f_m), x^*}
  = \int_{M^m} f_1(t_1)\cdot\ldots\cdot f_m(t_m)
     \ip{\phi(t_1, \ldots, t_m),x^*}
     \; d\mu^{\ot m}(t_1, \ldots, t_m).
 \end{align*}
We will not always notationally distinguish between a function
$\phi$ and the operator $T_\phi \in \g(L^2(M),E)$ that it
represents. The subspace of operators which can be represented by a
function is dense in $\g^m(L^2(M),E).$

\section{Wiener-It\^{o} chaos in Banach spaces} \label{sec:chaos}

In this section we will prove a Banach space analogue of the
classical Wiener-It\^o isometry. First we fix some notations.

Let $(\O,\F,\P)$ be a (sufficiently rich) probability space, let
 $H$ be a real separable Hilbert space and let $W: H\to L^2(\O)$ be an
isonormal Gaussian process on $H,$ i.e. $W$ maps $H$ isometrically
onto a closed subspace consisting of Gaussian random variables. We
assume that $\F$ is the $\s$-field generated by $\{W(h):h \in H\}.$
We fix an orthonormal basis $(u_j)_{j\geq1}$ of $H,$ and consider
the Gaussian sequence defined by $\g_j := W(u_j)$ for $j\geq 1.$
For $m \geq 0$ we consider the $m$-th Wiener-It\^o chaos,
 $$\H_m := \overline{\lin}\{ H_m(W(h)) : \|h\| =
1\},$$ where the closure is taken in $L^2(\O).$ Furthermore, let
$H^{\sot m}$ be the $m$-fold symmetric tensor power which is defined
to be the range of the orthogonal projection $P_\sot \in \L(H^{\hot
m})$ given by
 \begin{align*}
 P_\sot(h_1 \ot \cdots \ot h_m) =
   \frac{1}{m!} \sum_{\pi \in S_m} h_{\pi(1)} \ot \cdots \ot
   h_{\pi(m)}, \qquad h_1, \ldots, h_m \in H,
 \end{align*}
where $S_m$ is the group of permutations of $\{1,\ldots,m\}.$

A classical result of Wiener states that the following orthogonal
decomposition holds:
 \begin{align*}
L^2(\O,\F,\P) = \bigoplus_{m\geq 0} \H_m.
 \end{align*}
Moreover, the mapping $\Phi_m$ defined by
 \begin{align} \label{eq:canonical}
 \Phi_m :  P_{\sot}(u_{i_1} \ot \cdots \ot u_{i_m} )
       \mapsto (\i!/m!)^{1/2} \Psi_\i,
 \end{align}
extends to an isometry from $H^{\sot m}$ onto $\H_m.$  Recall that
$\Psi_\i$ is the generalised Hermite polynomial defined in Section
\ref{subsec:decoupling}, to which we refer for notations.

Let us  consider the vector-valued Gaussian chaos
 \begin{align*}
 \H_m(E) := \overline{\lin} \{ f \ot x : f \in \H_m,
                              x \in E \},
 \end{align*}
where the closure is taken in ${L^2(\O;E)}.$ The following
well-known result is a consequence of the decoupling result in
Theorem \ref{thm:decouplingD}(1) and the Kahane-Khintchine
inequalities. Extensive information on this topic can be found in
the monographs \cite{KW92,dlPG99}.

 \bprop \label{prop:KKchaos}
Let $E$ be a Banach space, let $m \geq 1,$ and let $1\leq
p,q<\infty.$ For all $F \in \H_m(E)$ we have $$\|F\|_{L^p(\O;E)}
   \eqsim_{m,p,q} \|F\|_{L^q(\O;E)}.$$
 \eprop

Our next goal is the construction of the spaces $\g^{\sot m}(H,E),$
which will be the Banach space substitutes for the symmetric Hilbert
space tensor powers.
We refer to Section \ref{subsec:gammanorms} for the definition of
the space $\g^m(H,E)$.
For $T \in \g^m(H,E)$  we define its symmetrisation $P_\sot T \in
\g^m(H,E)$ by
 \begin{align*}
 (P_\sot T ) (h_1, \ldots, h_m) :=
 \frac{1}{m!} \sum_{\pi\in S_m} T(h_{\pi(1)}, \ldots,
 h_{\pi(m)}), \quad h_1, \ldots, h_m \in H,
 \end{align*}
and we will say that $T \in \g^m(H,E)$ is symmetric if $P_\sot T =
T.$ The mapping $P_\sot$ is easily seen to be a projection in
$\L(\g^m(H,E))$ and we define $ \g^{\sot m }(H,E)$ to be its range.

We remark that if $K$ is a Hilbert space, then $\g^{\sot m}(H,K)$ is
isometrically isomorphic to the space $H^{\sot m} \hot K,$ where
$\hot$ denotes the Hilbert space tensor product.

Now we are ready to state the main result of this section, which is
a Banach space valued extension of the canonical isometry
\eqref{eq:canonical}.

 \bthm \label{thm:WienerIto}
Let $E$ be a Banach space, let $1\leq p<\infty,$ and let $m \geq 1.$
The mapping
  $$(\Phi_m \ot I) : P_\sot(h_{i_1} \ot \cdots \ot h_{i_m} )\ot x
       \mapsto (\i!/m!)^{1/2} \Psi_\i \ot x,$$
extends to a bounded operator $(\Phi_m \ot I) :  \g^{\sot m}(H,E)
\to L^p(\O;E),$ which maps $\g^{\sot m}(H,E)$ onto $\H_m(E).$
Moreover, we have equivalence of norms
 $$ \| (\Phi_m \ot I) T\|_{L^p(\O;E)} \eqsim_{m,p} \|T\|_{\g^{m}(H,E)},
 \quad T \in \g^{\sot m}(H,E).$$
 \ethm

 \bpr
Let $T$ be a symmetric operator of the form
\eqref{eq:elementaryoperator} and observe that
 $$ T = \sum_{|\i| = m, |\i|_\infty \leq n}
      P_{\sot}(u_{i_1} \ot \cdots \ot u_{i_m})\ot x_\i.$$
Using \eqref{eq:elementaryoperatornormgm}, the decoupling result
from Theorem \ref{thm:decouplingD}(1) and the Kahane-Khintchine
inequalities we obtain
 \begin{align*}
 \E \|(\Phi_m \ot I) T\|_{E}^p
 & = \E \Big\| \sum_{|\i|=m, |\i|_\infty \leq n}
 (\i!/m!)^{1/2} \Psi_\i    x_{\i}  \Big\|_E^p
 \\& \eqsim_{m,p} \E \Big\| \sum_{|\i| = m,
         |\i|_\infty \leq n}
 \g_{i_1}^{(1)} \cdot \ldots\cdot \g_{i_m}^{(m)}
   x_{\i}  \Big\|_E^p
 \eqsim_{m,p} \|T\|_{\g^m(H,E)}^p.
 \end{align*}
In view of Proposition \ref{prop:KKchaos} it is clear that $\Phi_m
\ot I$ maps $\g^{\sot m}(H,E)$ into $\H_m(E).$ To show that its
range is $\H_m(E),$ we observe that $\Phi_m \ot I (h^{\ot m} \ot x)
= H_m(W(h))\cdot x$ for all $h \in H$ with $\|h\|=1$ and all $x \in
E.$ Now the result follows from the norm estimate above and the
identity
 \begin{align*}
  \H_m(E) =  \overline{\lin} \{H_m(W(h)) \cdot x : \|h\| = 1,
   x \in E\},
 \end{align*}
where the closure is taken in $L^p(\O;E).$
 \epr

\brem In the special case that $E = \R$ and $p=2$ we recover the
classical Wiener-It\^o isometry. \erem

\brem \label{rem:Kconvex} Let $m \geq 1$ and let $J_m$ be the
orthogonal projection onto $\H_m.$ It is well known that for all
$1<p<\infty$ the restriction of $J_m$ to $L^p(\O) \cap L^2(\O)$
extends to a bounded projection on $L^p(\O).$  A Banach space $E$ is
said to be $K$-convex if $J_1 \ot I$ extends to a bounded operator on
$L^2(\O;E).$ Actually, this notion is usually defined using
Rademacher instead of Gaussian random variables, but this does not
affect the class of Banach spaces under consideration \cite{FTJ79}.
It has been shown by Pisier \cite{Pi82} that in this case the
operators $J_m \ot I$ (which will be denoted by $J_m$ below) are
bounded for all $m \geq 1$ and all $1<p<\infty.$ Every UMD space is
$K$-convex. These facts will be used in Sections \ref{sec:derivative}
and \ref{sec:OU}. \erem

\section{Multiple Wiener-It\^{o} integrals in Banach spaces}
\label{subsec:whitenoise} As in the previous section we consider a
real separable Hilbert space $H$ and an isometry $W: H\to L^2(\O)$
onto a closed linear subspace consisting of Gaussian random
variables.

In addition we assume in this section that $H = L^2(M,\B,\mu)$ for
some $\s$-finite non-atomic measure space $M.$ We let $\B_0 := \{ B
\in \B : \mu(B) < \infty\}.$ For $A \in \B_0$ we write with some
abuse of notation $W(A):= W(\one_A).$ In this way $W$ defines an
$L^2(\O)$-valued measure on $\B_0$ which is called the white noise
based on $\mu.$

Our next goal is to construct multiple stochastic integrals for
Banach space-valued functions. Our construction generalises the well
known multiple stochastic integral for Hilbert space valued
functions, and in another direction, the (single) stochastic
integral for Banach space valued functions which has been
constructed in \cite{vNW05}.

For fixed $m \geq 1$  we define $\eE_m(E)$ to be the linear space of
tetrahedral simple functions $F : M^m \to \R$ of the form
 \begin{align} \label{eq:elementaryfunction}
F = \sum_{|\i| = m, |\i|_\infty \leq n} \one_{A_{i_1} \times \cdots
\times A_{i_m}} \cdot x_{\i},
 \end{align}
where the $A_j$'s are pairwise disjoint sets in $\B_0,$ $n \geq 1,$
and the coefficients $x_{\i} \in E$ vanish whenever $j(\i) > 1$ for
some $j\geq 1.$
It is easy to see that such a function $F$ represents an operator
$T_F \in \g^m(L^2(M),E)$ in the sense described in Section
\ref{subsec:gammanorms}, and by taking an orthonormal basis
$(u_j)_{j \geq 1}$ of $L^2(M)$ with $u_j =
\mu(A_j)^{-1/2}\one_{A_j}$ for $j = 1, \ldots, n,$ one can check
that
 \beq\bal \label{eq:gammamnorm}
&\|T_F\|_{\g^m(L^2(M),E)}^2%
=
  \E  \Big\| \sum_{|\i|=m,|\i|_\infty \leq n} \g_{i_1}^{(1)}  \cdot \ldots \cdot
 \g_{i_m}^{(m)} \cdot \mu(A_1)^{1/2}  \cdot \ldots \cdot
\mu(A_n)^{1/2} \cdot x_\i \Big\|_E^2.
 \eal\eeq
We recall that $(\g_j^{(k)})_{j \geq 1}$ are independent Gaussian
sequences for $k \geq 1$.

 \blem \label{lem:tetraoperdense}
The collection of operators represented by functions in $\eE_m(E)$
is dense in $\g^m(L^2(M),E)$ for all $m \geq 1.$
 \elem

 \bpr This follows by reasoning as in the proof of the
corresponding scalar-valued result  \cite[p.10]{Nu06}, taking into
account that the measure space $M$ is non-atomic.
\epr

Suppose that $T_F \in \g^m(L^2(M),E)$ is represented by a strongly
measurable weakly-$L^2$ function $F.$ Then $T_F$ belongs to
$\g^{\sot m}(L^2(M),E)$ if and only if $F$ agrees almost everywhere
with its symmetrisation $\widetilde F$ defined by
 \begin{align*}
 \widetilde F (t_1, \ldots, t_m) :=
 \frac{1}{m!} \sum_{\pi\in S_m} F(t_{\pi(1)}, \ldots,
 t_{\pi(m)}).
 \end{align*}
For $F \in \eE_m(E)$ of the form \eqref{eq:elementaryfunction} we
define the multiple Wiener-It\^{o} integral $I_m(F) \in L^2(\O;E)$
by
 \begin{align} \label{eq:elementaryMWIintegral}
I_m(F) = \sum_ {|\i| = m,|\i|_\infty \leq n} W(A_{i_1})  \cdot
\ldots \cdot W( A_{i_m} ) \cdot x_\i.
 \end{align}
One easily checks that this definition does not depend on the
representation of $F$ as an element of $\eE_m(E).$ Moreover, $I_m$
is linear and $I_m(F) = I_m(\widetilde F).$  The next theorem may be
considered as a generalisation of the classical
Wiener-It\^{o}-isometry for multiple stochastic integrals to the
Banach space setting.

 \bthm \label{thm:multstochint}
Let $m \geq 1$ and $1\leq p <\infty.$ The operator $I_m : \eE_m(E)
\to L^p(\O;E)$ extends uniquely to a bounded operator $$I_m :
\g^m(L^2(M),E) \to L^p(\O;E),$$ which maps $\g^m(L^2(M),E)$ onto
$\H_m(E).$ Moreover, for all $F \in \g^m(L^2(M),E)$ we have:
 \beni
 \item $I_m F = I_m \widetilde F;$
 \item $\|I_m F\|_{L^p(\O;E)}
        \eqsim_{m,p} \| \widetilde F \|_{\g^m(L^2(M),E)}
        \leq \| F \|_{\g^m(L^2(M),E)}.$
 \een
 \ethm

 \bpr
First we show that for all $F \in \eE_m(E)$ the following
equivalence of norms holds:
 \begin{align*}
\|I_m F\|_{L^p(\O;E)} \eqsim_{m,p} \|\widetilde
F\|_{\g^m(L^2(M),E)}.
 \end{align*}
For that purpose we take  $F\in \eE_m(E)$ of the form
\eqref{eq:elementaryfunction}. Since $I_m(F) = I_m(\widetilde F)$ we
may assume that $F$ is symmetric, hence $ x_{(i_{\pi(1)},\ldots,
i_{\pi(m)}) } = x_{(i_1, \ldots, i_m) } $ for all permutations $\pi
\in S_m.$ Let $(u_j)_{j \geq 1}$ be an orthonormal basis of $L^2(M)$
with $u_j = \mu(A_j)^{-1/2}\one_{A_j}$ for $j = 1, \ldots , n,$ and
let $(\g_j)_{j \geq 1}$ be the Gaussian sequence $\g_j =W(u_j)$ for
$j \geq 1.$ Using  the decoupling inequalities from Theorem
\ref{thm:decouplingD}(2), \eqref{eq:gammamnorm}, and the
Kahane-Khintchine inequalities we obtain
\begin{align*}
\|I_m F\|_{L^p(\O,E)}^p
 & = \E\Big\|\sum_{|\i|=m,|\i|_\infty \leq n} W(A_{i_1})  \cdot \ldots
\cdot W( A_{i_m} ) \cdot x_\i \Big\|_E^p
 \\& =
  \E \Big\| \sum_{|\i|=m,|\i|_\infty \leq n} \g_{i_1}  \cdot \ldots
\cdot \g_{i_m} \cdot \mu(A_1)^{1/2}  \cdot \ldots \cdot
\mu(A_n)^{1/2} \cdot x_\i \Big\|_E^p
 \\& \eqsim_{m,p}
  \E  \Big\| \sum_{|\i|=m,|\i|_\infty \leq n}
    \g_{i_1}^{(1)}  \cdot \ldots
\cdot \g_{i_m}^{(m)} \cdot \mu(A_1)^{1/2}  \cdot \ldots \cdot
\mu(A_n)^{1/2} \cdot x_\i \Big\|_E^p
 \\& \eqsim_{m,p} \|F\|_{\g^m(L^2(M),E)}^p.
 \end{align*}
Now the first claim follows from Lemma \ref{lem:tetraoperdense}. To
prove that $I_m T \in \H_m(E)$ for all $T \in \g^m(L^2(M),E)$ we
first let $T = T_F$ for some tetrahedral function $F$  of the form
\eqref{eq:elementaryfunction}.
It follows from
\eqref{eq:elementaryMWIintegral} and the fact that
$$W(A_{j_1}) \cdot \ldots\cdot W(A_{j_m}) \in \H_m$$ whenever all
$j_k$'s are different, that $I_m T \in \H_m(E).$ Since $I_m$ is
continuous the same holds for general $T \in \g^m(L^2(M),E)$ by Lemma
\ref{lem:tetraoperdense}. To show that the mapping $I_m :
\g^m(L^2(M),E) \to \H_m(E)$ is surjective we proceed as in Theorem
\ref{thm:WienerIto}. The other statements are clear in view of Lemma
\ref{lem:tetraoperdense}.
 \epr

\section{The Malliavin derivative}\label{sec:derivative}

In this section we consider a complete probability space
$(\O,\F,\P),$ a real separable Hilbert space  $H,$ and an isonormal
Gaussian process $W: H \to L^2(\O).$ As before we assume that $\F$ is
the $\s$-algebra generated by $W.$

Let us introduce some notation. For $n\geq 1$ we denote by
$C_{pol}^{\infty}(\R^n)$ the vector space of all
$C^{\infty}$-functions $f:\R^n\to\R$ such that $f$ and its partial
derivatives of all orders have polynomial growth, i.e. for every
multi-index $\a$ there exist positive constants $C_{\a}, p_{\a}$
such that
$$|\partial_{\a} f(x) | \leq C_{\a}(1+|x|)^{p_\a}.$$
Let $\calS$ be the collection of all random variables $f:\O \to \R$
of the form
  \begin{align}\label{eq:cylindrical}
f = \varphi(W(h_1), \ldots, W(h_n))
  \end{align}
 for some $\varphi\in C_{pol}^{\infty}(\R^n)$, $h_1,\ldots,h_n\in H$
  and $n\geq1$.

For  a real Banach space  $E$ we consider the dense subspace
$\calS(E)$ of $L^p(\O;E),$ $1\leq p <\infty,$ consisting of all
functions $F:\O \to E$ of the form
$$F = \sum_{i=1}^n f_i \cdot x_i,$$ where $f_i \in \calS$ and $x_i \in E,$
$i=1,\ldots n.$
Occasionally it will be convenient to work with the space $\calP(E),$
which is defined similarly, except that the functions $\varphi$ are
required to be polynomials.

For a function $F = f \cdot x\in\calS(E)$ with $f$ of the form
\eqref{eq:cylindrical} we define its
 Malliavin derivative $DF$ by
 \begin{align}\label{eq:derivative}
 DF = \sum_{j=1}^n \partial_j \varphi (W(h_1), \ldots,
 W(h_n)) h_j \ot x.
 \end{align}
This definition extends to
$\calS(E)$ by linearity.
For $F \in \calS(E)$ the Malliavin derivative $DF$ is a random
variable which takes values in the algebraic tensor product $H \ot
E,$ which we endow with the norm $\| \cdot \|_{\g(H,E)}$ (cf.
Section \ref{subsec:gammanorms}).

The following result is the simplest case of the integration by parts
formula. We omit the proof, which is the same as in the scalar-valued
case \cite[Lemma 1.2.1]{Nu06}.

 \blem \label{lem:expDFh}
 If $F\in\calS(E)$, then
 $\E(DF(h)) = \E(W(h)F)$ for all $h\in H$.
 \elem

 A straightforward computation shows that the following product rule
holds:
 \begin{align*}
  D\ip{F,G} = \ip{DF,G} + \ip{F,DG},
  \quad F \in \calS(E), G \in \calS(E^*).
 \end{align*}
Here $\ip{\cdot,\cdot}$ denotes the duality between $E$ and $E^*.$
Combining this with Lemma \ref{lem:expDFh} we obtain the following
integration by parts formula:
 \begin{align}\label{eq:ibpprel}
 \E \ip{DF(h),G} = \E(W(h)\ip{F,G}) - \E\ip{F,DG(h)},
  \quad F \in \calS(E), G \in \calS(E^*).
 \end{align}

This identity is the main ingredient in the proof of the following
result which can be found in \cite[Proposition 2.3]{MvNCO}.

 \bprop \label{prop:closable}
The Malliavin derivative $D$ is closable as an operator from
$L^p(\O;E)$ into $L^p(\O;\g(H,E))$ for all $1 \leq p < \infty.$
 \eprop

With a slight abuse of notation we will denote the closure of $D$
again by $D.$ Its domain in $L^p(\O;E)$ will be denoted by
$\dD^{1,p}(\O;E),$ which is a Banach space endowed with the norm
 \begin{align*}
 \|F\|_{\dD^{1,p}(\O;E)} := \big( \|F\|_{L^p(\O;E)}^p
           + \|DF\|_{L^p(\O;\g(H,E))}^p  \big)^{1/p}.
 \end{align*}
Furthermore we will write $\dD^{1,p}(\O) := \dD^{1,p}(\O;\R).$

Derivatives of higher order are defined inductively. For  $n \geq 1$
we define
 \begin{align*}
  \dD^{n+1,p}(\O;E) :=& \{ F\in \dD^{n,p}(\O;E) : D^n F \in
                       \dD^{1,p}(\O;\g^n(H,E))\},\\
   D^{n+1}F      :=& D(D^n F), \quad  F\in \dD^{n+1,p}(\O;E).
 \end{align*}
It follows from Proposition \ref{prop:closable} that $D^n$ is a
closed and densely defined operator from $\dD^{n-1,p}(\O;E)$ into
$L^p(\O;\g^n(H,E)).$ Its domain is denoted by $\dD^{n,p}(\O;E)$
which is a Banach space endowed with the norm
 \begin{align*}
 \|F\|_{n,p} := \|F\|_{\dD^{n,p}(\O;E)} := \Big( \|F\|_{L^p(\O;E)}^p
           + \sum_{k=1}^n\|D^k F\|_{L^p(\O;\g^k(H,E))}^p  \Big)^{1/p}.
 \end{align*}
The main result in this section describes the behaviour of the
Malliavin derivative on the $E$-valued Wiener-It\^o chaoses. It
extends \cite[Proposition 1.2.2]{Nu06} to Banach spaces (and to
$1\leq p<\infty$, but this is well-known in the scalar case).

 \bthm \label{thm:Dchaosm}
Let $E$ be a Banach space, let $1\leq p<\infty$ and let $m\geq 1.$
Then we have $\H_m(E) \subset \dD^{1,p}(\O;E)$ and $D(\H_m(E))
\subset \H_{m-1}(\g(H,E)).$ Moreover, the following equivalence of
norms holds:
 \begin{align*}
\|DF\|_{L^p(\O;\g(H,E))} \eqsim_{p,m} \|F\|_{L^p(\O;E)}, \quad F \in
\H_m(E).
 \end{align*}
 \ethm

 \bpr
Let $(u_j)_{j\geq 1}$ be an orthonormal basis of $H$ and put $\g_j
:= W(u_j).$ Let $(\g_j^{(k)})_{j \geq 1}$ and $(\widetilde
\g_j)_{j\geq1}$ be independent copies of $(\g_j)_{j \geq 1}.$ For
$\i = (i_1, \ldots, i_m)$ and $k \geq 1$ we will write $(\i,k) =
(i_1, \ldots, i_m,k).$

First we take $F \in \H_m(E)$ of the form
$$F = \sum_{ |\i| = m, |\i|_\infty \leq n } \tfrac{\i!}{{m!}^{1/2}}
 \prod_{j=1}^n H_{j(\i)}(\g_j)\,x_\i.$$
  Clearly we may assume without loss of generality
that the coefficients $x_\i$ are symmetric,
 i.e. $x_\i = x_{\i'}$ whenever $\i'$ is a permutation of $\i.$

It follows from Theorem \ref{thm:decouplingD}(1) that
 \beq\bal\label{eq:Fchaosm}
 \E\|F\|_E^p
  &= \E \bigg\| \sum_{|\i| = m,|\i|_\infty \leq n}\tfrac{\i!}{{m!}^{1/2}}
      \prod_{j=1}^n H_{j(\i)}(\g_j) x_\i  \bigg\|_E^p
  \\&\eqsim_{m,p} \E
   \bigg\| \sum_{|\i| = m}
   \g_{i_1}^{(1)} \cdot \ldots \cdot  \g_{i_m}^{(m)}  x_\i
   \bigg\|_E^p.
 \eal\eeq
On the other hand, by a change of variables to modify the range of
summation from $\{ |\i| = m\}$ to $\{ |\i| = m-1\},$ and rearranging
terms, we obtain with the convention that $H_{-1}=0,$
 \begin{align*}
 DF &=  \sum_{ |\i| = m,|\i|_\infty \leq n } \tfrac{\i!}{{m!}^{1/2}}
\sum_{k=1}^n
    \prod_{j\neq k} H_{j(\i)}(\g_j)
     H_{k(\i) - 1}(\g_k) \cdot u_k  \ot x_\i,
   \\& =  \sum_{k=1}^n u_k \ot \bigg(  m^{1/2}
     \sum_{ |\i| = {m-1},|\i|_\infty \leq n } \tfrac{\i!}{(m-1)!}
      \prod_{j= 1}^n  H_{j(\i)}(\g_j) x_{(\i,k)}\bigg)
   \\& =  \sum_{k=1}^n u_k \ot \bigg( \tfrac{m^{1/2}}{(m-1)!^{1/2}}
     \sum_{ |\i| = {m-1},|\i|_\infty \leq n } \tfrac{\i!^{1/2}}{(m-1)!^{1/2}}
      \Psi_\i  x_{(\i,k)}\bigg).
 \end{align*}
Using the Kahane-Khintchine inequalities and Theorem
\ref{thm:decouplingD}(1) once more, we find
 \beq\bal\label{eq:DFchaosm}
 \E\|DF\|_{\g(H,E)}^p
  & \eqsim_p \E  \widetilde\E \Big\| \sum_{k=1}^n
      \widetilde\g_k DF(u_k)
          \Big\|_E^p
 \\ & = \E \widetilde\E \Big\| \tfrac{m^{1/2}}{(m-1)!^{1/2}}
      \sum_{k=1}^n \widetilde \g_k
     \sum_{ |\i| = m-1,|\i|_\infty \leq n }  \tfrac{\i!^{1/2}}{(m-1)!^{1/2}}
      \Psi_\i  x_{(\i,k)} \Big\|_E^p
 \\& = \widetilde\E \E \Big\| \tfrac{m^{1/2}}{(m-1)!^{1/2}}
     \sum_{ |\i| = m-1,|\i|_\infty \leq n }
     \tfrac{\i!^{1/2}}{(m-1)!^{1/2}}
      \Psi_\i
 \Big(\sum_{k=1}^n \widetilde \g_k x_{(\i,k)} \Big)\Big\|_E^p
 \\&  \eqsim_{m,p} \widetilde\E \E
   \Big\|  \sum_{ |\i| = m-1,|\i|_\infty \leq n } \sum_{k=1}^n
     \g_{i_1}^{(1)}\cdot \ldots \cdot \g_{i_{m-1}}^{(m-1)}
 \widetilde \g_k x_{(\i,k)} \Big\|_E^p.
   \eal\eeq
Comparing \eqref{eq:Fchaosm} and \eqref{eq:DFchaosm} yields the norm
estimate. The theorem follows by the closedness of $D$ and the fact
that functions $F$ of the form considered above are dense in
$\H_m(E)$.
 \epr

 \brem
In the special case that $E$ is a UMD Banach space the result above
is known. Indeed, it follows from Meyer's inequalities (Theorem
\ref{thm:Meyerineq}) that
$$\|D F\|_{L^p(\O;\g(H,E))} \eqsim_{p,E} m^{1/2}\|F\|_{L^p(\O;E)}, \quad F \in \H_m(E).$$
This formula gives an explicit dependence on $m,$ but in contrast
with Theorem 5.3 the constants depend on (the Hilbert transform
constants of) $E.$ We return to this observation in Section
\ref{sec:OU}.
 \erem

\section{Meyer's inequalities and its consequences} \label{sec:OU}

Let $(P(t))_{t \geq 0} \subset \L(L^2(\O))$ be the Ornstein-Uhlenbeck semigroup defined by
 \begin{align}
 P(t) := \sum_{m\geq 0} e^{-mt} J_m.
 \end{align}
As is well known, this semigroup extends to a $C_0$-semigroup of positive contractions on $L^p(\O)$ for all $1\leq p < \infty.$ We refer the reader to \cite{Nu06} for proofs of these and other elementary properties.

Let  $E$ be an arbitrary Banach space. By positivity of $P$, $(P(t)
\ot I)_{t \geq 0}$ extends to a $C_0$-semigroup of contractions on
the Lebesgue-Bochner spaces $L^p(\O;E)$ for $1 \leq p < \infty$ which
will be denoted by $(P_E(t))_{t \geq 0}.$ The domain in $L^p(\O;E)$
of its infinitesimal generator $L_E$ is denoted $\Dom_p(L_E).$ The
subordinated semigroup $(Q_E(t))_{t \geq 0}$ is defined by
 \begin{align} \label{eq:subordinated}
 Q_E(t) f := \int_0^\infty P_E(s) f \;d\nu_t(s),
 \end{align}
where the probability measure $\nu_t$ is given by
 \begin{align}
d\nu_t(s) =  \frac{t}{2\sqrt{\pi s^3}}e^{-t^2/4s}\, ds, \quad t > 0.
 \end{align}
The generator of $(Q_E(t))_{t \geq 0}$ will be denoted by $C_E.$ As
is well known we have $$C_E = -(-L_E)^{1/2}.$$ Often, when there is
no danger of confusion, we will omit the subscripts $E.$

The next lemma is a vector-valued analogue of the representation of
$L$ as a generator associated with a Dirichlet form. We omit the
proof which follows from the scalar-valued analogue in a
straightforward way.

 \blem \label{lem:dirichletform}
Let $E$ be a UMD space. For all $F \in \calP(E)$ and $G \in
\dD^{1,p}(\O;E^*)$ we have
 \begin{align*}
 \E\ip{L_E F,G} = \E[DF,DG]_{\g}.
 \end{align*}
 \elem

In the following Lemma we collect some useful commutation relations,
which follow easily from the corresponding scalar-valued results.

 \blem \label{lem:commutation}
Let $E$ be a Banach space and let $1 \leq p < \infty.$ \beni
\item For $F \in \dD^{1,p}(\O;E)$ we have $P_E(t) F \in \dD^{1,p}(\O;E)$
      and $$D P_E(t) f = e^{-t}P_{\g(H,E)} DF.$$
\item For $F \in \dD^{1,p}(\O;E)$ we have $Q_E(t) F \in \dD^{1,p}(\O;E)$
      and $$D Q_E(t) f = Q_{\g(H,E)}^{(1)} DF,$$
      where $Q_{\g(H,E)}^{(1)}$ is the semigroup  generated by
      $-(I-L_{\g(H,E)})^{1/2}.$
\item For $F \in \calP(E)$ we have $L_E F \in \dD^{1,p}(\O;E)$ and
      $D L_E F = -(I-L_{\g(H,E)}) DF.$
\item For $F \in \calP(E)$ we have $C_E F \in \dD^{1,p}(\O;E)$ and
      $D C_E F = -(I - L_{\g(H,E)})^{1/2} DF.$
\een
 \elem

Pisier proved in \cite{Pi88} that Meyer's inequalities extend to UMD
spaces. Formulated in the language of $\gamma$-norms his result reads
as follows.

 \bthm[Meyer's inequalities] \label{thm:Meyerineq}
Let $E$ be a UMD  space and let $1 < p < \infty.$ Then $\Dom_p(C_E) =
\dD^{1,p}(\O;E)$ and for all $f\in \dD^{1,p}(\O;E)$ the following
two-sided estimate holds:
 \begin{align} \label{eq:Meyerineq}
  \| C_E f \|_{L^p(\O;E)}
   \eqsim_{p,E} \|D f\|_{L^p(\O;\g(H,E))}.
 \end{align}
 \ethm

In Theorem \ref{thm:meyergeneral} we shall state an extension of this
result.

The following lemma is the crucial ingredient in the proof of Meyer's
multiplier Theorem. The proof in the scalar case in \cite[Lemma
1.4.1]{Nu06} does not extend to the vector-valued setting, since it
depends heavily on the Hilbert space structure of $L^2(\O).$ We give
a simple proof in the case that $E$ is a UMD space, which is based on
Meyer's inequalities. Recall that $J_m$ denotes the chaos projection
considered in Remark \ref{rem:Kconvex}.

 \blem \label{lem:crucial}
Let $1<p<\infty$ and let $E$ be a UMD space. For each $N \geq 1$ and
$t > 0$ we have
 \begin{align*}
 \| P(t) (I - J_0 - J_1 - \ldots - J_{N-1}) \|_{\L(L^p(\O;E))}
  \lesssim_{E,p,N} e^{-Nt}.
 \end{align*}
 \elem

 \bpr
For $F \in \calP(E)$ we set
 \begin{align*}
  R F = D \sum_{m=1}^\infty m^{-1/2} J_m F, \qquad
  S \Big( D \sum_{m=0}^\infty J_m F \Big)
      := \sum_{m=1}^\infty m^{1/2} J_m F.
 \end{align*}
Note that the sums consists of finitely many terms since $F \in
\calP(E).$ Both operators are well-defined and bounded by Theorem
\ref{thm:Meyerineq}.
 Using the fact that
 \begin{align*}
S^N R^N F = \sum_{m=N}^\infty J_m F,
 \end{align*}
we obtain by Lemma \ref{lem:commutation} and Theorem
\ref{thm:Meyerineq},
 \begin{align*}
 \| P(t) &(I - J_0 - J_1 - \ldots - J_{N-1}) F \|_{L^p(\O;E)}
 \\&   = \Big\|\sum_{m=N}^\infty e^{-mt} J_m F \Big\|_{L^p(\O;E)}
  = \|S^{N}R^{N} P(t) F \|_{L^p(\O;E)}
 \\& =  \|S^{N} e^{-Nt}P(t) R^{N} F \|_{L^p(\O;E)}
  \leq   e^{-Nt} \|S\|^N \| R\|^{N} \| F \|_{L^p(\O;E)}.
  \end{align*}
 \epr

Using this lemma, the remainder of the proof of Meyer's multiplier
Theorem \cite{Mey84} in the scalar case as given in \cite[Theorem
1.4.2]{Nu06} extends verbatim to the vector-valued setting. It is
even possible to allow operator-valued multipliers.

 \bthm[Meyer's Multiplier Theorem] \label{thm:multiplier}
Let $1<p<\infty,$ let $E$ be a UMD space,  and let
$(a_k)_{k=0}^\infty \subset \L(L^p(\O;E))$ be a sequence of bounded
linear operators such that $\sum_{k=0}^\infty
\|a_k\|_{\L(L^p(\O;E))} N^{-k} < \infty$ for some $N \geq 1.$ If
$(\phi(n))_{n\geq 0}\subset \L(L^p(\O;E))$ is a sequence of
operators satisfying $\phi(n) := \sum_{k=0}^\infty a_k n^{-k}$ for
$n \geq N,$ then the operator $T_\phi$ defined by
 $$T_\phi F  := \sum_{n=0}^\infty \phi(n) J_n F,
     \quad F \in \calP(E)$$
extends to a bounded operator on $L^p(\O;E).$
 \ethm

 As a first application of the multiplier theorem we determine the
spectrum of $L.$ We start with a simple but useful lemma.

 \blem \label{lem:Jmseparates}
Let $E$ be a $K$-convex Banach space, let $1<p<\infty,$ and let $F
\in L^p(\O;E)$ such that $J_m F = 0$ for all $m \geq 0.$ Then $F =
0$ in $L^p(\O;E).$
 \elem

 \bpr
For $G \in \calP(E^*)$ we have $$\E\ip{F, G} = \E\ip{F, \sum_{m \geq
0} J_m G} =  \E\ip{ \sum_{m \geq 0} J_m F, G} = 0.$$ This implies
the result, since $\calP(E^*)$ is dense in $L^q(\O;E^*),$ hence
weak$^*$-dense in $L^p(\O;E)^*.$

 \epr

 \bprop \label{prop:spectrum}
Let $1< p< \infty$ and let $E$ be a UMD space. Then
 \begin{align*}
 \s(-L) = \{0,1,2, \ldots\}.
 \end{align*}
Moreover, every integer $m \geq 0$ is an eigenvalue of $-L$ and
$\ker (m + L) = \H_m(E).$
 \eprop

 \bpr
To prove that $\{0,1,2,\ldots\} \subset \s(-L)$ we take an integer $m
\geq 0$ and a non-zero $F \in \H_m(E).$ Since $P(t) F = e^{-mt} F$ it
follows that $ F\in \Dom_p(L)$ and $(m+L) F = 0,$ hence $m \in
\s(-L)$ and $\ker (m + L) \supset \H_m(E).$

To show the converse inclusion for the spectrum, take $\l \in \C
\setminus \{0,1,2,\ldots\}.$
To prove that $\l + L$ is injective, take $F \in \ker(\l+L).$ Since
$J_m$ is bounded for $m \geq 0$ by Remark \ref{rem:Kconvex} (UMD
spaces are $K$-convex), it follows that $J_m LF = LJ_m F = -m J_m
F.$ This implies that $(\l - m )J_m F = J_m(\l + L)F = 0,$ hence
$J_m F = 0$ for all $m \geq 0,$ so that $F = 0$ by Lemma
\ref{lem:Jmseparates}.

To prove surjectivity, we conclude from the Multiplier Theorem
\ref{thm:multiplier} that
$$ R_{\l} := \sum_{m=0}^\infty \frac{1}{\l-m} J_m$$
extends to a bounded operator on $L^p(\O;E).$ Using the fact that $L$
is closed, we infer that $(\l + L)R_{\l} = I,$ hence $\l+L$ is
surjective.

It remains to show that $\ker (m + L) \subset \H_m(E)$ for all
 $m\geq 0.$ Take $F \in \ker (m + L).$ Since
$$(m-k)J_k F = (m+L)J_k F = J_k(m+L)F = 0$$ for all integers
$k \geq 0,$ we have $J_k F = 0$ for all $k \neq m.$ This implies
that $J_k(F - J_m F) = 0,$ hence $F = J_m F \in \H_m(E)$ by Lemma
\ref{lem:Jmseparates}.

 \epr

Next we give the general form of Meyer's inequalities in the
language of $\g$-radonifying norms.
This result is stated in a slightly different setting in
\cite[Theorem 1.17]{MN93}, but the proof given there contains
a gap. More precisely, the last formula for the function
$\psi$ defined in \cite[p.300]{MN93} should be replaced by
$\psi(t) = \frac12
e^{-t/2}(I_0(\frac{t}{2})+I_1(\frac{t}{2})).$ This function
however is not contained in $L^1(0,\infty);$ but this is
needed to conclude the proof.

The proof given below uses Lemma \ref{lem:crucial}, which is based
on the first order Meyer inequalities from Theorem
\ref{thm:Meyerineq}. This allows us to adapt the argument in the
scalar case from \cite[Theorem 1.5.1]{Nu06}.
%

 \bthm[Meyer's inequalities, general case] \label{thm:meyergeneral}
Let $E$ be a UMD space, let $1 < p < \infty$ and let $n \geq 1.$ Then
$\Dom_p(C^{n}) = \dD^{n,p}(\O;E),$ and  for all $F
\in\dD^{n,p}(\O;E)$ we have
 \beq\bal \label{eq:Meyergeneral}
  \|D^n F\|_{L^p(\O;\g^n(H,E))}
     &\lesssim_{p,E,n} \|C^n F \|_{L^p(\O;E)}
     \\&\lesssim_{p,E,n} \|F\|_{L^p(\O;E)}
       + \|D^n F\|_{L^p(\O;\g^n(H,E))}.
 \eal\eeq
 \ethm

 \bpr
The proof proceeds by induction. The case $n=1$ has  been treated in
Theorem \ref{thm:Meyerineq}. Suppose that \eqref{eq:Meyergeneral}
holds for some $n \geq 1.$ Using Lemma \ref{lem:commutation} and the
fact that the operator $C^n(I-L)^{-n/2} = (-L)^{n/2}(I-L)^{-n/2}$ is
bounded on $L^p(\O;E)$ we obtain by the induction hypothesis
 \begin{align*}
            \E\|   D^{n+1} F     \|_{\g^{n+1}(H,E)}^p
 & \lesssim_{p,E,n} \E\|      C^{n} D F  \|_{\g(H,E)}^p
   \lesssim_{p,E,n} \E\|  (I-L)^{n/2} D F\|_{\g(H,E)}^p
\\&=        \E\|  D C^{n}       F\|_{\g(H,E)}^p
   \eqsim_{p,E}   \E\|  C^{n+1}       F\|_{E}^p.
 \end{align*}

To prove the second inequality, we note that according to Remark
\ref{rem:Kconvex},
 $$\|C^n(J_0 + \ldots +J_{n-1})F\|_p
     \lesssim_{p,E,n} \| F \|_p, \quad F \in L^p(\O;E).$$
Therefore it suffices to show by induction that $$\|C^n F
\|_{L^p(\O;E)} \lesssim_{p,E,n} \|D^n F\|_{L^p(\O;\g^n(H,E))}$$ for
all $F \in \calP(E)$ with $J_0 F= \ldots = J_{n-1} F = 0.$

Let us assume that this statement holds for some $n\geq1$ and take
$F \in \calP(E)$ satisfying $J_0 F = \ldots = J_{n}F = 0.$ It
follows from Lemma \ref{lem:crucial} that $(P(t))_{t\geq0}$
restricts to a $C_0$-semigroup $(P_n(t))_{t\geq0}$ on $$X_{n,p}(E)
:= \overline{\bigoplus_{m\geq n} \H_m(E)}^{L^p(\O;E)},
$$
satisfying the growth bound  $\|P_n(t)\|_{\L(X_{n,p}(E))}
\lesssim_{E,p,n} e^{-nt}$ for some constant $K$ depending on $n.$
Consequently (see e.g. \cite[Proposition 3.8.2]{Ar01}), we have
$$\| (\a-L)^{1/2} F \|_p \eqsim_{p,E} \| (\b-L)^{1/2} F \|_p, \quad F \in
X_{n,p}(E),$$ for all $\a,\b>-n,$ and in particular is
$(I-L)^{1/2}C^{-1}$ bounded on $X_{n,p}(E).$ Using Lemma
\ref{lem:commutation} and the fact that $C^n DF \in
X_{n,p}(\g(H,E)),$ it follows that
 \begin{align*}
 \E\|C^{n+1}F\|_E^{p}
  &\eqsim_{p,E}  \E\|DC^{n}F\|_{\g(H,E)}^{p}
  =       \E\|(I-L)^{n/2}DF\|_{\g(H,E)}^{p}
  \\&\lesssim_{p,E,n} \|(I-L)^{n/2}C^{-n}\|_{\L(X_{n,p}(\g(H,E)))}^p
      \E\|C^n DF\|_{\g(H,E)}^{p}
  \\&\eqsim_{p,E,n} \E\|D^{n+1}F\|_{\g^{n+1}(H,E)}^{p}.
 \end{align*}
 \epr

As an application of Meyer's inequalities we will show that
$\g(H,E)$-valued Malliavin differentiable random variables are
contained in the domain of the divergence operator $\d.$ First we
give the precise definition of $\d.$

Fix an exponent $1< p<\infty$ and let $\tfrac{1}{p}+\tfrac{1}{q}=1.$
For the moment let $D$ denote the Malliavin derivative on
$L^q(\O;E^*)$, which is a densely defined closed operator with domain
$\dD^{1,q}(\O;E^*)$ and taking values in $L^q(\O;\g(H,E^*))$. We let
the domain $\Dom_p(\d)$ consist of all $u \in L^p(\O;\g(H,E))$ for
which there exists an $F_u \in L^p(\O;E)$ such that  $$\E[u,DG]_\g =
\E\ip{F_u,G}
 \  \text{ for all } G \in \dD^{1,q}(\O;E^*).$$
The function $F_u$, if it exists, is uniquely determined. We set
$$\d(u):= F_u , \quad X \in \Dom_p(\d).$$
In other words,  $\d$ is the part of the adjoint operator $D^*$ in
$L^p(\O;\g(H,E))$ which maps into $L^p(\O;E).$ Here we identify
 $L^p(\O;\g(H,E))$ and $L^p(\O;E)$ in a natural way with
subspaces of $(L^q(\O;\g(H,E^*)))^*$ and $(L^q(\O;E^*))^*$
respectively.

The divergence operator $\d$ is easily seen to be closed and densely
defined. For more information we refer to \cite{MvNCO}. The proof of
the following result is a  variation of the proof of the
scalar-valued result in \cite[Proposition 1.5.4]{Nu06}.

 \bprop \label{prop:skorohodbounded}
Let $1<p<\infty$ and let $E$ be a UMD space. The operator $\d$ is
bounded  from $\dD^{1,p}(\O;\g(H,E))$ into $L^p(\O;E).$
 \eprop

 \bpr
Let $u \in \dD^{1,p}(\O;\g(H,E))$ and $G \in \calP(E^*).$
Using Theorem \ref{thm:Dchaosm} we find that $\|D J_1 G\|_p \eqsim_p
\|J_1 G\|_p,$ and therefore
 \beq\bal \label{eq:divone}
  \E[ u, D(J_0+J_1)G]_\g
  &\leq \|u\|_{L^p(\O;\g(H,E))} \|D(J_0+J_1)G\|_{L^q(\O;\g(H,E^*))}
  \\&\lesssim_{p,E} \|u\|_{L^p(\O;\g(H,E))}  \|G\|_{L^q(\O;E^*)}.
 \eal\eeq
Now we assume that $J_0 G = J_1 G = 0.$ By the Multiplier Theorem
\ref{thm:multiplier} the operator
 \begin{align*}
 T := \sum_{m=2}^\infty \frac{m}{m-1} J_m
 \end{align*}
is bounded on $L^p(\O;\g(H,E)).$ By Lemma \ref{lem:crucial} the
operator $L^{-1}$ is well defined on $X_{1,p}(E),$ where we use the
notation from the proof of Theorem \ref{thm:meyergeneral}. This
justifies the use of $L^{-1}$ in the following computation. Using
Lemma \ref{lem:dirichletform} and Theorem \ref{thm:meyergeneral} we
obtain
 \beq\bal \label{eq:divtwo}
    \E [ u, DG]_{\g}
 &=  \E [ u, LL^{-1}DG]_{\g}
 =  \E [D u, DL^{-1}DG]_{\g}
 \\&\leq \|D u\|_{L^p(\O;\g^2(H,E))}
      \| DL^{-1}DG\|_{L^q(\O;\g^2(H,E^*))}
 \\&= \|D u\|_{L^p(\O;\g^2(H,E))}
      \| D^2 L^{-1} T G\|_{L^q(\O;\g^2(H,E^*))}
 \\&\lesssim_{p,E} \|D u\|_{L^p(\O;\g^2(H,E))}
      \|  G\|_{L^q(\O;E^*)}.
 \eal\eeq
Combining \eqref{eq:divone} and \eqref{eq:divtwo} we conclude that
for all $G \in \calP(E^*)$ we have
$$\E[ u, DG]_{\g}
  \lesssim_{p,E} \|u\|_{\dD^{1,p}(\O;E)}  \|G\|_{L^q(\O;E^*)}.
 $$
It follows that there exists an $F_u \in (L^q(\O;E^*))^*$ such that
$\E[ u, DG]_{\g} = \E\ip{F_u,G}.$ Since $E$ is a UMD space, we
conclude that $F_u \in L^p(\O;E)$ and we obtain the desired result.
 \epr

 For $1\leq p < \infty$ we define the vector space of exact
$E$-valued processes as
$$ L_e^p(\O;\g(H,E)) = \{DF : F \in \dD^{1,p}(\O;E)\}. $$
The next result is concerned with the representation of random
variables as divergences of exact processes.

 \bprop
Let $E$ be a UMD space, let $1<p<\infty$ and let $F \in L^p(\O;E).$
Then $U := D L^{-1}(F - \E(F)) $ is the unique element in
$L_e^p(\O;\g(H,E))$ satisfying $$F = \E(F) + \d(U).$$
 \eprop

 \bpr
By an easy computation we see that
 \begin{align}\label{eq:Fformula}
F = \E(F) + \d D(L^{-1}(F - \E(F)))
 \end{align}
for all $F\in \calP(E).$ It follows from Lemma \ref{lem:crucial} (or
Proposition \ref{prop:spectrum}) that $L^{-1}$ is well-defined and
bounded on $\{G \in L^p(\O;E) : \E(G) = 0 \}.$  Meyer's inequalities
imply that $D$ is bounded from $\Dom_p(L)$ into
$\dD^{1,p}(\O;\g(H,E)),$ and by Proposition
\ref{prop:skorohodbounded} we have that $\d$ is bounded from
$\dD^{1,p}(\O;\g(H,E))$ into $L^p(\O;E).$ Using these facts and an
approximation argument with elements from $\calP(E)$ we conclude
that the right hand side of \eqref{eq:Fformula} is well-defined for
all $F \in L^p(\O;E),$ and the identity remains valid.

To prove uniqueness, suppose that $F = \E(F) + \d(DF')$ for some $F'
\in \dD^{1,p}(\O;E)$ with $DF' \in \Dom_p(\d),$ and put $G := F' -
L^{-1}(F - \E(F)).$ Then $\d D G =  0,$ hence $\ip{ G , LP}=0$ for
all polynomials $P \in \calP(E^*).$ In particular, for all $m \geq
1$ and all $P \in \calP(E^*) \cap \H_m(E^*)$ one has $ \ip{ G ,
mP}=0,$ and since $\calP(E^*) \cap \H_m(E^*)$ is dense in
$\H_m(E^*),$ we have $\ip{J_m G, \widetilde F} = \ip{G, J_m
\widetilde F} = 0$ for all $\widetilde F \in L^q(\O;E^*).$ It
follows that $J_m G = 0$ for all $m \geq 1,$ which implies $J_m (G -
J_0 G) = 0$ for all $m \geq 0.$ We conclude that $G = J_0 G$ by
Lemma \ref{lem:Jmseparates}, hence $F' = L^{-1}(F - \E F) + x$ for
some $x \in E.$ We conclude that $DF' = DL^{-1}(F - \E F),$ which is
the desired identity.
 \epr

We conclude the paper with an application of the vector-valued
Malliavin calculus developed in this work. We give a new proof of
Theorem \ref{thm:decouplingD}(1) under the additional assumption that
$E$ is a UMD space, which is based on Meyer's inequalities. This
approach seems to be new even in the scalar-valued case.

 \bthm \label{thm:decouplingUMD}
Let $E$ be a UMD space, let $1<p<\infty,$ and define $F$ and
$\widetilde F$ as in Theorem \ref{thm:decouplingD}(1). Then we have
 \begin{align*}
\|F\|_p \eqsim_{p,m,E} \|\widetilde F\|_p.
 \end{align*}
 \ethm

 \bpr
We argue as in the proof of Theorem \ref{thm:Dchaosm}. By
\eqref{eq:Fchaosm} we have
  $$ \E\|F\|_E^p
  = \E \Big\| \sum_{|\i| = m,|\i|_\infty \leq n}\tfrac{\i!^{1/2}}{m!^{1/2}}
   \Psi_\i x_\i   \Big\|_E^p
 $$
and according to \eqref{eq:DFchaosm},
 $$ \E\|DF\|_{\g(H,E)}^p
 \eqsim_p \widetilde\E \E \Big\| \tfrac{m^{1/2}}{(m-1)!^{1/2}}
   \sum_{ |\i| = m-1,|\i|_\infty \leq n }
    \tfrac{\i!^{1/2}}{(m-1)!^{1/2}}
      \Psi_\i
 \Big(\sum_{k=1}^n \widetilde \g_k x_{(\i,k)} \Big)\Big\|_E^p.
   $$
Noting that  $CF = m^{1/2}F,$ Meyer's inequalities imply that
\begin{align*}\E \Big\| &\sum_{|\i| = m,|\i|_\infty \leq
n}\tfrac{\i!^{1/2}}{m!^{1/2}}
  \Psi_\i  x_\i  \Big\|_E^p
 \\& \eqsim_{p,m,E}
  \widetilde\E \E \Big\|  \sum_{ |\i| = m-1,|\i|_\infty \leq n }
    \tfrac{\i!^{1/2}}{(m-1)!^{1/2}}
      \Psi_\i
 \Big(\sum_{k=1}^n \widetilde \g_k x_{(\i,k)} \Big)\Big\|_E^p.
   \end{align*}
The desired result is obtained by repeating this procedure $m-1$
times.
 \epr

\bibliographystyle{ams-pln}
\bibliography{UMDMalliavin}

\def\lfhook#1{\setbox0=\hbox{#1}{\ooalign{\hidewidth
  \lower1.5ex\hbox{'}\hidewidth\crcr\unhbox0}}}
\providecommand{\bysame}{\leavevmode\hbox to3em{\hrulefill}\thinspace}
\begin{thebibliography}{10}

\bibitem{AG93}
M.~A. Arcones and E.~Gin{\'e}, \emph{On decoupling, series expansions, and tail
  behavior of chaos processes}, J. Theoret. Probab. \textbf{6} (1993), no.~1,
  101--122.

\bibitem{Ar01}
W.~Arendt, C.J.K. Batty, M.~Hieber, and F.~Neubrander, \emph{Vector-valued
  {L}aplace transforms and {C}auchy problems}, Monographs in Mathematics,
  vol.~96, Birkh\"auser Verlag, Basel, 2001.

\bibitem{Bo}
V.I. Bogachev, \emph{Gaussian measures}, Mathematical Surveys and Monographs,
  vol.~62, American Mathematical Society, Providence, RI, 1998.

\bibitem{CT06}
R.~A. Carmona and M.~R. Tehranchi, \emph{Interest rate models: an infinite
  dimensional stochastic analysis perspective}, Springer Finance,
  Springer-Verlag, Berlin, 2006.

\bibitem{dlPG99}
V.~H. de~la Pe{\~n}a and E.~Gin{\'e}, \emph{Decoupling}, Probability and its
  Applications (New York), Springer-Verlag, New York, 1999, From dependence to
  independence, Randomly stopped processes. $U$-statistics and processes.
  Martingales and beyond.

\bibitem{dlPMS95}
V.~H. de~la Pe{\~n}a and S.~J. Montgomery-Smith, \emph{Decoupling inequalities
  for the tail probabilities of multivariate {$U$}-statistics}, Ann. Probab.
  \textbf{23} (1995), no.~2, 806--816.

\bibitem{FTJ79}
T.~Figiel and N.~Tomczak-Jaegermann, \emph{Projections onto {H}ilbertian
  subspaces of {B}anach spaces}, Israel J. Math. \textbf{33} (1979), no.~2,
  155--171.

\bibitem{hytonen:lps}
T.~Hytonen, \emph{{Littlewood--Paley--Stein theory for semigroups in UMD
  spaces}}, Rev. Mat. Iberoamericana, to appear.

\bibitem{IW89}
N.~Ikeda and S.~Watanabe, \emph{Stochastic differential equations and diffusion
  processes}, second ed., North-Holland Mathematical Library, vol.~24,
  North-Holland Publishing Co., Amsterdam, 1989.

\bibitem{KWsf}
N.J. Kalton and L.~Weis, \emph{The ${H}^\infty$-functional calculus and square
  function estimates}, in preparation.

\bibitem{KW04}
P.C. Kunstmann and L.~Weis, \emph{Maximal {$L\sb p$}-regularity for parabolic
  equations, {F}ourier multiplier theorems and {$H\sp \infty$}-functional
  calculus}, Functional analytic methods for evolution equations, Lecture Notes
  in Math., vol. 1855, Springer, Berlin, 2004, pp.~65--311.

\bibitem{Kw87}
S.~Kwapie{\'n}, \emph{Decoupling inequalities for polynomial chaos}, Ann.
  Probab. \textbf{15} (1987), no.~3, 1062--1071.

\bibitem{KW92}
S.~Kwapie\'{n} and W.~A. Woyczy\'{n}ski, \emph{Random series and stochastic
  integrals: single and multiple}, Probability and its Applications,
  Birkh\"auser Boston Inc., Boston, MA, 1992.

\bibitem{MvNCO}
J.~Maas and J.M.A.M.~van Neerven, \emph{A {C}lark-{O}cone formula in {UMD}
  {B}anach spaces}, arXiv: 0709.2021.

\bibitem{Ma78}
P.~Malliavin, \emph{Stochastic calculus of variation and hypoelliptic
  operators}, Proceedings of the International Symposium on Stochastic
  Differential Equations (Res. Inst. Math. Sci., Kyoto Univ., Kyoto, 1976) (New
  York), Wiley, 1978, pp.~195--263.

\bibitem{MN93}
P.~Malliavin and D.~Nualart, \emph{Quasi-sure analysis and {S}tratonovich
  anticipative stochastic differential equations}, Probab. Theory Related
  Fields \textbf{96} (1993), no.~1, 45--55.

\bibitem{MWZ05}
E.~Mayer-Wolf and M.~Zakai, \emph{The divergence of {B}anach space valued
  random variables on {W}iener space}, Probab. Theory Related Fields
  \textbf{132} (2005), no.~2, 291--320, Correction available at
  arXiv:0710.4483.

\bibitem{McCT86}
T.~R. McConnell and M.~S. Taqqu, \emph{Decoupling inequalities for multilinear
  forms in independent symmetric random variables}, Ann. Probab. \textbf{14}
  (1986), no.~3, 943--954.

\bibitem{McCT87}
\bysame, \emph{Decoupling of {B}anach-valued multilinear forms in independent
  symmetric {B}anach-valued random variables}, Probab. Theory Related Fields
  \textbf{75} (1987), no.~4, 499--507.

\bibitem{Mey84}
P.-A. Meyer, \emph{Transformations de {R}iesz pour les lois gaussiennes},
  Seminar on probability, XVIII, Lecture Notes in Math., vol. 1059, Springer,
  Berlin, 1984, pp.~179--193.

\bibitem{vNVW07}
J.M.A.M.~van Neerven, M.C. Veraar, and L.~Weis, \emph{{Stochastic integration
  of processes with values in a UMD Banach space}}.

\bibitem{vNW07}
J.M.A.M.~van Neerven and L.~Weis, \emph{Stochastic integration of
  operator-valued functions with respect to {B}anach space-valued {B}rownian
  motion}, preprint.

\bibitem{vNW05}
\bysame, \emph{Stochastic integration of functions with values in a {B}anach
  space}, Studia Math. \textbf{166} (2005), no.~2, 131--170.

\bibitem{Nu06}
D.~Nualart, \emph{The {M}alliavin calculus and related topics}, second ed.,
  Probability and its Applications (New York), Springer-Verlag, Berlin, 2006.

\bibitem{Pi78}
G.~Pisier, \emph{Some results on {B}anach spaces without local unconditional
  structure}, Compositio Math. \textbf{37} (1978), no.~1, 3--19.

\bibitem{Pi82}
\bysame, \emph{Holomorphic semigroups and the geometry of {B}anach spaces},
  Ann. of Math. (2) \textbf{115} (1982), no.~2, 375--392.

\bibitem{Pi88}
\bysame, \emph{Riesz transforms: a simpler analytic proof of {P}.-{A}.
  {M}eyer's inequality}, S\'eminaire de Probabilit\'es, XXII, Lecture Notes in
  Math., vol. 1321, Springer, Berlin, 1988, pp.~485--501.

\bibitem{Pis89}
\bysame, \emph{The volume of convex bodies and {B}anach space geometry},
  Cambridge Tracts in Mathematics, vol.~94, Cambridge University Press,
  Cambridge, 1989.

\bibitem{Sh94}
I.~Shigekawa, \emph{Sobolev spaces of {B}anach-valued functions associated with
  a {M}arkov process}, Probab. Theory Related Fields \textbf{99} (1994), no.~3,
  425--441.

\end{thebibliography}

\end{document}